\documentclass[11pt]{amsart}

\usepackage{amsmath,amsthm}
\usepackage{amssymb}
\usepackage{euscript}

\numberwithin{equation}{section} 

\newcommand{\sqsp}{\renewcommand{\baselinestretch}{1.1}\tiny\normalsize}

\newtheorem{thm}[equation]{Theorem}
\newtheorem{cor}[equation]{Corollary}
\newtheorem{lemma}[equation]{Lemma}
\newtheorem{prop}[equation]{Proposition}

\theoremstyle{definition}
\newtheorem{remark}[equation]{Remark}

\newcommand{\bC}{\mathbf{C}}

\newcommand{\bP}{\mathbf{P}}
\newcommand{\bH}{\mathbf{H}}
\newcommand{\bZ}{\mathbf{Z}}

\DeclareMathOperator{\holim}{holim}
\DeclareMathOperator{\Hom}{Hom}
\DeclareMathOperator{\genus}{Genus}
\DeclareMathOperator{\Aut}{Aut}
\DeclareMathOperator{\dg}{deg}           
\DeclareMathOperator{\LCM}{LCM}          

\newcommand{\cpinfty}{\bC\bP^\infty}     
\newcommand{\hpinfty}{\bH\bP^\infty}     

\begin{document}
\title[Genus of infinite quaternionic projective space]
{Maps to spaces in the genus of infinite quaternionic projective space}
\author[Donald Yau]{Donald Yau}

\begin{abstract}
Spaces in the genus of infinite quaternionic projective space which admit essential maps from infinite complex projective space are classified.  In these cases the sets of homotopy classes of maps are described explicitly.  These results strengthen the classical theorem of McGibbon and Rector on maximal torus admissibility for spaces in the genus of infinite quaternionic projective space.  An interpretation of these results in the context of Adams-Wilkerson embedding in integral $K$-theory is also given.
\end{abstract}



\date{\today}

\maketitle

\sqsp


\section{Introduction and statement of results}
\label{sec:intro}

In an attempt to understand Lie groups through their classifying spaces, Rector \cite{rector} classified the genus of $\hpinfty$, the infinite projective space over the quaternions, considered as a model for the classifying space $BS^3$.  The homotopy type of a spaces $X$ is said to be in the genus of $\hpinfty$, denoted $X \in \genus(\hpinfty)$, if the $p$-localizations of $X$ and $\hpinfty$ are homotopy equivalent for each prime $p$.  One often speaks of a space rather than its homotopy type when considering genus.  Rector's classification \cite{rector} of the genus of $\hpinfty$ is as follows.
\bigskip


\begin{thm}[Rector]
\label{rector}
Let $X$ be a space in the genus of $\hpinfty$.  Then for each prime $p$ there exists a homotopy invariant $(X/p) \in \lbrace \pm 1 \rbrace$ such that the following statements hold.
\begin{enumerate}
\item The $(X/p)$ for $p$ primes provide a complete list of homotopy classification invariants for the genus of $\hpinfty$.  
\item Any combination of values of the $(X/p)$ can occur.  In particular, the genus of $\hpinfty$ is uncountable.  
\item The invariant $(\hpinfty/p)$ is $1$ for all primes $p$.
\item The space $X$ has a maximal torus if and only if $X$ is homotopy equivalent to $\hpinfty$.
\end{enumerate}
\end{thm}
\bigskip

The invariant $(X/p)$ is now known as the \emph{Rector invariant at the prime} $p$.  Actually, for the last statement about the maximal torus, Rector only proved it for the odd primes.  That is, if $X$ has a maximal torus, then $(X/p)$ is equal to $1$ for all odd primes $p$.  Then McGibbon \cite{mcgibbon} proved it for the prime $2$ as well.  Here $X$ is said to have a maximal torus if there exists a map from $\cpinfty$, the infinite complex projective space, to $X$ whose homotopy theoretic fiber has the homotopy type of a finite complex.

For a space $X$ in the genus of $\hpinfty$ which is not homotopy equivalent to $\hpinfty$, the nonexistence of a maximal torus does not rule out the possibility that there could be some essential (that is, non-nullhomotopic) maps from $\cpinfty$ to $X$.  The main purposes of this paper are (1) to describe spaces in the genus of $\hpinfty$ for which this can happen (Theorem \ref{thm:null}) and (2) to compute the maps in these cases (Theorem \ref{thm:maps}).

The following is our first main result, which classifies spaces in the genus of $\hpinfty$ which admit essential maps from $\cpinfty$. 
\bigskip


\begin{thm}\label{thm:null}
Let $X$ be a space in the genus of $\hpinfty$.  Then the following statements are equivalent. 
\begin{enumerate}
\item There exists an essential map from $\cpinfty$ to $X$.
\item There exists a nonzero integer $k$ such that $(X/p) = (k/p)$ for all but finitely many primes $p$.
\item There exists a cofinite set of primes $L$ such that $\hpinfty$ and $X$ become homotopy equivalent after localization at $L$.
\end{enumerate}
\end{thm}
\bigskip

Here $(k/p)$ is the Legendre symbol of $k$, which is defined whenever $p$ does not divide $k$.  If $p$ is odd and if $p$ does not divide $k$, then 
$(k/p) = 1$ (resp.\ $-1$) if $k$ is a quadratic residue (resp.\ non-residue) mod $p$.  If $p = 2$ and if $k$ is odd, then 
$(k/2) = 1$ (resp.\ $-1$) if $k$ is a quadratic residue (resp.\ non-residue) mod $8$.

Before discussing related issues, let us first record the following immediate consequence of Theorem \ref{thm:null}.


\bigskip
\begin{cor}
\label{cor1:null}
There exist only countably many homotopically distinct spaces in the genus of $\hpinfty$ which admit essential maps from $\cpinfty$.
\end{cor}
\bigskip

Indeed, each nonzero integer $k$ can determine only countably many homotopically distinct spaces $X$ in the genus of $\hpinfty$ satisfying the second condition in Theorem \ref{thm:null}.

The second condition of Theorem \ref{thm:null} gives an arithmetic description of spaces in the genus of $\hpinfty$ which occur as the target of essential maps from $\cpinfty$.  Since it involves Rector invariants, it is specific to the genus of $\hpinfty$ and is not very convenient for generalizations.  The last condition of Theorem \ref{thm:null}, on the other hand, is geometric and is more suitable for possible generalizations of the theorem.

Having characterized spaces in the genus of $\hpinfty$ which admit nontrivial maps from $\cpinfty$, we proceed to compute the maps themselves.  Now for any space $X$ in the genus of $\hpinfty$, the $K$-theory $K(X)$ of $X$, as a filtered ring, is a power series ring $\bZ \lbrack \lbrack b^2 u_X \rbrack \rbrack$ (see Proposition \ref{prop1:null}), where $u_X$ is some element in $K^4(X)$ and $b$ is the Bott element in $K^{-2}(\text{pt})$.  Here, and throughout the rest of the paper, $K(-)$ denotes complex $K$-theory with coefficients over the integers $\bZ$.  So if $f \colon \cpinfty \to X$ is any map, then its induced map in $K$-theory defines an integer $\dg(f)$, called the \emph{degree of} $f$, by the equation
\begin{equation}
\label{eq:degree}
f^*(b^2 u_X) 
~=~ \dg(f)(b\xi)^2 + \text{higher terms in } b\xi,
\end{equation}
where $b\xi$ is the ring generator in the power series ring $K(\cpinfty) = \bZ \lbrack \lbrack b\xi \rbrack \rbrack$.  Note that $\dg(f)$ is, up to a sign, simply the degree of the induced map of $f$ in integral homology in dimension $4$.  According to a result of Dehon and Lannes \cite{dl}, the homotopy class of such a map $f$ is determined by its degree.  We will therefore identify such a map with its degree in the sequel.  The degrees of a self-map of $X$, a self-map of the $p$-localization $\hpinfty_{(p)}$ of $\hpinfty$, or a map from $\cpinfty$ to $\hpinfty_{(p)}$ can be defined similarly.

To describe the maps from $\cpinfty$ to $X \in \genus(\hpinfty)$ up to homotopy, we need only describe the possible degrees of such maps.  Let's first consider the classical case.  There is a maximal torus inclusion $i \colon \cpinfty \to \hpinfty$ of degree $1$, and any other map $f \colon \cpinfty \to \hpinfty$ factors through $i$ up to homotopy.  A special case of a theorem of Ishiguro, M$\o$ller, and Notbohm \cite[Thm.\ 1]{imn} says that for any space $X$ in the genus of $\hpinfty$, the degrees of essential self-maps of $X$ consist of precisely the squares of odd numbers.  For the classical case, $X = \hpinfty$, this result is due to Sullivan \cite[p.\ 58-59]{sullivan}.  Therefore, the degrees of essential maps from $\cpinfty$ to $\hpinfty$ also consist of precisely the odd squares.

The situation in general is quite similar.  Recall that any space $X$ in the genus of $\hpinfty$ can be obtained as a homotopy inverse limit \cite{bk} 
\begin{equation}
\label{eq:construction}
X ~=~ \holim_q \, \left\lbrace \hpinfty_{(q)} \buildrel r_q \over \longrightarrow \hpinfty_{(0)} \buildrel n_q \over \longrightarrow \hpinfty_{(0)}\right\rbrace.
\end{equation}
Here $q$ runs through all primes, $r_q$ is the natural map from the $q$-localization to the rationalization of $\hpinfty$, and $n_q$ is an integer relatively prime to $q$, satisfying $(n_q/q) = (X/q)$.  The integer $n_2$ also satisfies $n_2 \equiv 1$ (mod $4$).

Now if $X$ admits an essential map from $\cpinfty$, and thus satisfies the condition in Theorem \ref{thm:null} for some nonzero integer $k$, then the integers $n_q$ can be chosen so that the set $\lbrace n_q \colon q \text{ primes} \rbrace$ contains only finitely many distinct integers.  So it makes sense to talk about the least common multiple of the integers $n_q$, denoted $\LCM(n_q)$.  Now we define an integer $T_X$ as 
\begin{equation}
\label{eq:lcm}
T_X ~=~ \min \lbrace \LCM(n_q) \, \colon \, X = \holim_q\, (n_q \circ r_q)\rbrace .
\end{equation}
That is, choose the integers $n_q$ as in \eqref{eq:construction} so as to minimize their least common multiple, and $T_X$ is defined to be $\LCM(n_q)$.

We are now ready for the second main result of this paper, which describes the maps from $\cpinfty$ to $X \in \genus(\hpinfty)$.


\bigskip
\begin{thm}
\label{thm:maps}
Let $X$ be a space in the genus of $\hpinfty$ which admits an essential map from $\cpinfty$.  Then the following statements hold.
\begin{enumerate}
\item There exists a map $i_X \colon \cpinfty \to X$ of degree $T_X$.
\item The map $i_{\hpinfty} \colon \cpinfty \to \hpinfty$ is the maximal torus inclusion.
\item Given any map $f \colon \cpinfty \to X$, there exists a self-map $g$ of $X$ such that $f$ is homotopic to $g \circ i_X$.
\item The degrees of essential maps from $\cpinfty$ to $X$ are precisely the odd squares multiples of $T_X$.
\end{enumerate}
\end{thm}
\bigskip

It should be noted that the integer $T_X$ does \textbf{not} determine the homotopy type of $X$.  For example, consider the spaces $X$ and $Y$ in the genus of $\hpinfty$ with Rector invariants
\begin{equation}
(X/p) ~=~ \begin{cases}1 & \text{if } p \not= 3 \\
                      -1 & \text{if } p = 3\end{cases},
\qquad
(Y/p) ~=~ \begin{cases}1 & \text{if } p \not= 5 \\
                      -1 & \text{if } p = 5.\end{cases}
\end{equation}
Then, of course, $X$ is not homotopy equivalent to $Y$ because their Rector invariants at the prime $3$ are distinct.  But it is easy to see that $T_X = 2 = T_Y$.

Theorems \ref{thm:null} and \ref{thm:maps} are closely related to the (non)existence of Adams-Wilkerson type embeddings of finite $H$-spaces in integral $K$-theory.  As mentioned before, a map $f \colon \cpinfty \to X \in \genus(\hpinfty)$ is essential if and only if $\dg(f)$ is nonzero.  Thus, if there exists an essential map from $\cpinfty$ to $X \in \genus(\hpinfty)$, then $K(X)$ can be embedded into $K(\cpinfty)$ as a sub-$\lambda$-ring.  The converse is also true.  Indeed, a theorem of Notbohm and Smith \cite[Thm.\ 5.2]{ns} says that the function
\[
\alpha \colon \lbrack \cpinfty, X \rbrack ~\to~
\Hom_\lambda(K(X),K(\cpinfty))
\]
which sends (the homotopy class of) a map to its induced map in $K$-theory, is a bijection.  (Here $\lbrack -,- \rbrack$ and $\Hom_\lambda(-,-)$ denote, respectively, sets of homotopy classes of maps between spaces and of $\lambda$-ring homomorphisms.)  So a $\lambda$-ring embedding $K(X) \to K(\cpinfty)$ must be induced by an essential map from $\cpinfty$ to $X$.  Therefore, Theorem \ref{thm:null} and Corollary \ref{cor1:null} can be restated in this context as follows.


\bigskip
\begin{thm}
\label{thm':null}
Let $X$ be a space in the genus of $\hpinfty$.  Then $K(X)$ can be embedded into $K(\cpinfty)$ as a sub-$\lambda$-ring if, and only if, there exists a nonzero integer $k$ such that $(X/p) = (k/p)$ for all but finitely many primes $p$.  This is true if, and only if, there exists a cofinite set of primes $L$ such that $\hpinfty$ and $X$ become homotopy equivalent after localization at $L$.

In particular, there exist only countably many homotopically distinct spaces $X$ in the genus of $\hpinfty$ whose $K$-theory $\lambda$-rings can be embedded into that of $\cpinfty$ as a sub-$\lambda$-ring.
\end{thm}
\bigskip

Before Theorem \ref{thm':null}, there is at least one space in the genus of $\hpinfty$ whose $K$-theory $\lambda$-ring was known to be non-embedable into the $K$-theory $\lambda$-ring of $\cpinfty$.  This example was due to Adams \cite[p.\ 79]{adams}.

\bigskip
\begin{remark}
Theorems \ref{thm:null} and \ref{thm:maps} can also be regarded as an attempt to understand the set of homotopy classes of maps from $X$ to $Y$, where $X \in \genus(BG)$ and $Y \in \genus(BK)$ with $G$ and $K$ some connected compact Lie groups.  This problem, especially the case $G = K = S^3 \times \cdots \times S^3$, was studied extensively by Ishiguro, M$\o$ller, and Notbohm \cite{imn}. 
\end{remark}
\bigskip

This finishes the presentation of our main results.  The rest of this paper is organized as follows.  Theorems \ref{thm:null} and \ref{thm:maps} are proved in \S \ref{sec:proof} and \S \ref{sec:proof of thm maps}, respectively.   For use in \S \ref{sec:proof}, we recall in  \S \ref{sec:filtered ring} some results from \cite{adic}, with sketches of proofs, about $K$-theory filtered rings.


\section*{Acknowledgment}
The author would like to thank Dietrich Notbohm for his interest and comments in this work, and to the referee whose suggestions, especially the last condition of Theorem \ref{thm:null}, are very helpful.  The author would like to express his deepest gratitude to his thesis advisor Haynes Miller for many hours of stimulating conversations and encouragement.


\section{$K$-theory filtered ring}
\label{sec:filtered ring}

In preparation for the proof of Theorem \ref{thm:null} in \S \ref{sec:proof}, a result from \cite{adic} about $K$-theory filtered ring is reviewed in this section.

To begin with, a \emph{filtered ring} is a pair $(R, \{ I^R_n \})$ consisting of: 
\begin{enumerate}
\item A commutative ring $R$ with unit; 
\item A decreasing filtration $R = I^R_0 \supset I^R_1 \supset \cdots$ of ideals of $R$ such that $I^R_i I^R_j \subset I^R_{i+j}$ for all $i, j \geq 0$.  
\end{enumerate}
A map between two filtered rings is a ring homomorphism which preserves the filtrations.  With these maps as morphisms, the filtered rings form a category.

Every space $Z$ of the homotopy type of a CW complex gives rise naturally to an object $(K(Z), \{K_n(Z)\})$, which is usually abbreviated to $K(Z)$, in this category.  Here $K(Z)$ and $K_n(Z)$ denote, respectively, the complex $K$-theory ring of $Z$ and the kernel of the restriction map $K(Z) \to K(Z_{n-1})$, where $Z_{n-1}$ denotes the $(n-1)$-skeleton of $Z$.  Using a different CW structure of $Z$ will not change the filtered ring isomorphism type of $K(Z)$, as can be easily seen by using the cellular approximation theorem.  The symbol $K^r_s(Z)$ denotes the subgroup of $K^r(Z)$ consisting of elements whose restrictions to $K^r(Z_{s-1})$ are equal to $0$.

The following result, which is proved in \cite{adic}, will be needed in \S \ref{sec:proof}.  For the reader's convenience we include here a sketch of the proof.  In what follows $b \in K^{-2}(\text{pt})$ will denote the Bott element.

\bigskip
\begin{prop}
\label{prop1:null}
Let $X$ be a space in the genus of $\hpinfty$.  Then the following statements hold.
\begin{enumerate}
\item There exists an element $u_X \in K^4_4(X)$ such that $K(X) = \bZ\lbrack \lbrack b^2 u_X \rbrack \rbrack$ as a filtered ring.
\item For any odd prime $p$, the Adams operation $\psi^p$ satisfies
\begin{equation}\label{eq2:null}
\psi^p\left(b^2 u_X\right) 
~=~ \left(b^2 u_X\right)^p + 2 \, (X/p) \, p \left(b^2 u_X\right)^{(p+1)/2} + pw + p^2 z,
\end{equation}
where $w$ and $z$ are some elements in $K^0_{2p+3}(X)$ and $K^0_4(X)$, respectively.  In particular, we have
\begin{equation}\label{eq3:null}
\psi^p\left(b^2 u_X\right) 
~=~ 2\, (X/p)\, p \left(b^2 u_X\right)^{(p+1)/2} \quad \left(\text{mod } K^0_{2p+3}(X) \text{ and } p^2\right).
\end{equation}
\end{enumerate}
\end{prop}
\bigskip

\begin{proof}[Sketch of the proof of Proposition \ref{prop1:null}]
For the first assertion, Wilkerson's proof of the classification theorem \cite[Thm.\ I]{wilkerson} of spaces of the same $n$-type for all $n$ can be easily adapted to show the following.  There is a bijection between the following two pointed sets:
\begin{enumerate}
\item The pointed set of isomorphism classes of filtered rings $(R,\{I^R_n\})$ with the properties:
  \begin{enumerate} 
  \item The natural map $R \to \varprojlim_n\, R/I^R_n$ is an isomorphism, and 
  \item $R/I^R_n$ and $K(\hpinfty)/K_n(\hpinfty)$ are isomorphic as filtered rings for all $n > 0$.
  \end{enumerate}
\item The pointed set $\varprojlim^1_n\, \Aut(K(\hpinfty)/K_n(\hpinfty))$.
\end{enumerate}  
Here $\Aut(-)$ denotes the group of filtered ring automorphisms, and the $\varprojlim^1$ of a tower of not-necessarily abelian groups is as defined in \cite{bk}.  It is not difficult to check that the hypothesis on $X$ implies that the filtered ring $K(X)$ has the above two properties, (a) and (b).  Moreover, by analyzing the subquotients $K_n(\hpinfty)/K_{n+1}(\hpinfty)$, one can show that the map 
\[
\Aut(K(\hpinfty)/K_{n+1}(\hpinfty)) ~\to~ \Aut(K(\hpinfty)/K_n(\hpinfty))
\]
is surjective for each integer $n$ greater than $4$.  The point is that any automorphism of $K(\hpinfty)/K_n(\hpinfty)$ can be lifted to an endomorphism of $K(\hpinfty)/K_{n+1}(\hpinfty)$ without any difficulty.  Then, since the quotient $K(\hpinfty)/K_{n+1}(\hpinfty)$ is a finitely generated abelian group, one only has to observe that the chosen lift is surjective.  Thus, the above $\varprojlim^1$ is the one-point set, and hence $K(X)$ and $K(\hpinfty)$ are isomorphic as filtered rings.  This establishes the first assertion.

The second assertion concerning the Adams operations $\psi^p$ is an easy consequence of the first assertion, Atiyah's theorem \cite[Prop.\ 5.6]{atiyah}, and the definition of the Rector invariants $(X/p)$ \cite{rector}.
\end{proof}


\section{Proof of Theorem \ref{thm:null}}
\label{sec:proof}
In this section the proof of Theorem \ref{thm:null} is given.

Recall that the complex $K$-theory of $\cpinfty$ as a filtered $\lambda$-ring is given by $K(\cpinfty) = \bZ \lbrack \lbrack b\xi \rbrack \rbrack$ for some $\xi \in K^2_2(\cpinfty)$, where $b \in K^{-2}(\text{pt})$ is the Bott element.  The Adams operations on the generator are given by
\begin{equation}\label{eq1:null}
\psi^r(b\xi) ~=~ (1 + b\xi)^r - 1 \quad (r = 1, 2, \ldots).
\end{equation}

Fix a space $X$ in the genus of $\hpinfty$ and write $\bZ \lbrack \lbrack b^2 u_X \rbrack \rbrack$ for its $K$-theory filtered ring (cf.\ Proposition \ref{prop1:null}). 

We will prove Theorem \ref{thm:null} by proving the implications (1) $\Rightarrow$ (2) $\Rightarrow$ (3) $\Rightarrow$ (1).  Each implication is contained in one subsection below.


\subsection{Proof of (1) implies (2)}
\label{subsec:proof of (1) implies (2)}

This part of Theorem \ref{thm:null} is contained in the next Lemma.


\bigskip
\begin{lemma}\label{lem1:null}
Let $p$ be an odd prime and $k$ be a nonzero integer relatively prime to $p$.  If there exists an essential map $f \colon \cpinfty \to X$ of degree $k$, then $(X/p) = (k/p)$.
\end{lemma}
\bigskip

\begin{proof}
We will compare the coefficients of $(b\xi)^{p+1}$ in the equation
\begin{equation}\label{eq4:null}
f^* \psi^p\left(b^2 u_X\right) 
~=~ \psi^p f^*\left(b^2 u_X\right) \quad \left(\text{mod } K^0_{2p+3}(\cpinfty) \text{ and } p^2\right).
\end{equation}
Working modulo $K^0_{2p+3}(\cpinfty)$ and $p^2$, it follows from \eqref{eq:degree} and \eqref{eq3:null} that
\[
\begin{split}
f^* \psi^p\left(b^2 u_X\right) 
& =  2 \, (X/p) \, p \left(k b^2 \xi^2\right)^{(p+1)/2} \\
& =  2\, (X/p)\, p\, k^{(p+1)/2} (b\xi)^{p+1}.
\end{split}
\]
Similarly, still working modulo $K^0_{2p+3}(\cpinfty)$ and $p^2$, it follows from \eqref{eq:degree} and \eqref{eq1:null} that
\[
\begin{split}
\psi^p f^*\left(b^2 u_X\right) 
& =  k \psi^p(b^2\xi^2) \\ 
& =  k \psi^p(b\xi)^2 \\ 
& =  2pk (b\xi)^{p+1}.
\end{split}
\]
Thus, we obtain the congruence relation
\begin{equation}
\label{eq6:null}
2\, (X/p)\, p\, k^{(p+1)/2} 
~\equiv~ 2pk \quad (\text{mod }p^2).
\end{equation}
Since $(k/p)$ is congruent to $ k^{(p-1)/2}$ (mod $p$) (see, for example, \cite[Thm.\ 3.12]{nathanson}) and since $p$ is odd and relatively prime to $k$, \eqref{eq6:null} is equivalent to the congruence relation
\begin{equation}
(X/p) \, (k/p)  
~\equiv~ 1 \quad (\text{mod } p).
\end{equation}
Hence $(X/p) = (k/p)$, as desired.

This finishes the proof of Lemma \ref{lem1:null}.
\end{proof}

This shows that (1) implies (2) in Theorem \ref{thm:null}.
\bigskip

\subsection{Proof of (2) implies (3)}
\label{subsec:proof of (2) implies (3)}

Suppose that there exists a nonzero integer $k$ such that $(X/p) = (k/p)$ for all primes $p$, except possibly $p_1, \ldots, p_s$.  The prime factors of $k$ are among the $p_i$.  Let $L$ be the set consisting of all primes except the $p_i$, $1 \leq i \leq s$.  We will show that $\hpinfty$ and $X$ become homotopy equivalent after localization at $L$.

First note that for any space $Y$ in the genus of $\hpinfty$ and for any subset $I$ of primes, the $I$-localization of $Y$ can be obtained as
\begin{equation}
\label{eq:L-localization of Y}
Y_{(I)} ~=~ \holim_{q \in I}\, \left\lbrace \hpinfty_{(q)} \xrightarrow{n_q \, \circ \, r_q} \hpinfty_{(0)}\right\rbrace.
\end{equation}
In particular, we have
\begin{equation}
\label{eq:X_L}
X_{(L)} ~=~ \holim_{q \in L}\, \left\lbrace \hpinfty_{(q)} \xrightarrow{k \, \circ \, r_q} \hpinfty_{(0)}\right\rbrace
\end{equation}
and
\begin{equation}
\label{eq:BSU(2)_L}
\hpinfty_{(L)} ~=~ \holim_{q \in L}\, \left\lbrace \hpinfty_{(q)} \xrightarrow{r_q} \hpinfty_{(0)}\right\rbrace.
\end{equation}
Now for each prime $q \in L$, let $f_q$ be a self-map of $\hpinfty_{(q)}$ of degree $k^{-1}$.  Since $k$ is a $q$-local unit (because $q$ does not divide $k$), it is easy to see that each $f_q$ is a homotopy equivalence.  Moreover, the two maps
\begin{equation}
r_q,\, k \circ r_q \circ f_q \colon \hpinfty_{(q)} \to \hpinfty_{(0)}
\end{equation}
coincide.  Therefore, the maps $f_q$ $(q \in L)$ glue together to yield a map
\begin{equation}
f \colon \hpinfty_{(L)} \to X_{(L)}
\end{equation}
which is a homotopy equivalence, since each $f_q$ is.

This shows that (2) implies (3) in Theorem \ref{thm:null}.

\bigskip


\subsection{Proof of (3) implies (1)}
\label{subsec:proof of (3) implies (1)}

Suppose that there exists a cofinite set of primes $L$ such that $\hpinfty_{(L)}$ and $X_{(L)}$ are homotopy equivalent.  Write $p_1, \ldots, p_s$ for the primes not in $L$, and write $r_L$ for the natural map from $X_{(L)}$ to $\hpinfty_{(0)}$.

To construct an essential map from $\cpinfty$ to $X$, first note that $X$ can be constructed as the homotopy inverse limit of the diagram
\begin{equation}
  \hpinfty_{(p_i)} \xrightarrow{r_{p_i}} 
  \hpinfty_{(0)} \xrightarrow{n_{p_i}} \hpinfty_{(0)}
  \xleftarrow{r_L} X_{(L)}
\end{equation}
in which $i$ runs from $1$ to $s$.

For each $i$, $1 \leq i \leq s$, let $f_{p_i}$ be a map from $\cpinfty$ to $\hpinfty_{(p_i)}$ of degree $M/n_{p_i}$, where $M = \prod_{i=1}^s\, n_{p_i}$.  Also, let $f_L$ denote a map from $\cpinfty$ to $X_{(L)}$ of degree $M$, which exists because $X_{(L)}$ has the same homotopy type as $\hpinfty_{(L)}$.  It is then easy to see that the two maps
\begin{equation}
r_L \circ f_L, \, n_{p_i} \circ r_{p_i} \circ f_{p_i} \colon \cpinfty \to \hpinfty_{(0)}
\end{equation}  
coincide for any $1 \leq i \leq s$.  Therefore, the maps $f_{p_i}$ $(1 \leq i \leq s)$ and $f_L$ glue together to yield an essential map 
\begin{equation}
f \colon \cpinfty \to X
\end{equation}
through which all the maps $f_{p_i}$ and $f_L$ factor.

This shows that (3) implies (1) in Theorem \ref{thm:null}.

The proof of Theorem \ref{thm:null} is complete.

\bigskip


\section{Proof of Theorem \ref{thm:maps}}
\label{sec:proof of thm maps}
 
In this section we prove Theorem \ref{thm:maps}.

Fix a space $X$ in the genus of $\hpinfty$ which admits an essential map from $\cpinfty$.

First we note that part (2) follows from the discussion preceding Theorem \ref{thm:maps}, since it is obvious that the integer $T_{\hpinfty}$ is $1$.

Since any essential self-map of $X$ is a rational equivalence, part (4) is an immediate consequence of parts (1) and (3) and a result of Ishiguro, M$\o$ller, and Notbohm \cite[Thm.\ 1]{imn} which says that the degrees of essential self-maps of $X$ are precisely the odd squares.

Now we consider part (1).  Suppose that the integers $n_q$ as in \eqref{eq:construction} are chosen so that there are only finitely many distinct integers in the set $\lbrace n_q \colon q \text{ primes}\rbrace$ and that $T_X$ is their least common multiple (see \eqref{eq:lcm} for the definition of $T_X$).  Denote by $l_q \colon \hpinfty \to \hpinfty_{(q)}$ the $q$-localization map and by $i \colon \cpinfty \to \hpinfty$ the maximal torus inclusion.  A self-map of $\cpinfty$ of degree $m$ on $H^2(\cpinfty;\bZ)$ is simply denoted by $m$.  Now for each prime $q$ define a map $f_q \colon \cpinfty \to \hpinfty_{(q)}$ to be the composition
\begin{equation}
\label{eq:f_q}
\cpinfty \xrightarrow{M/n_q} 
\cpinfty \buildrel i \over \to 
\hpinfty \xrightarrow{l_q} \hpinfty_{(q)}.
\end{equation}
It is then easy to see that the two maps
\begin{equation}
n_q \circ r_q \circ f_q, \, n_{q^\prime} \circ r_{q^\prime} \circ f_{q^\prime} \colon \cpinfty \to \hpinfty_{(0)}
\end{equation}
coincide for any two primes $q$ and $q^\prime$.  Therefore, the maps $f_q$ glue together to yield an essential map
\begin{equation}
\label{eq:the map f}
f \colon \cpinfty \to X
\end{equation}
through which every map $f_q$ factors.  The map $f$ has degree $T_X$ because its induced map in rational cohomology in dimension $4$ does.

Finally, for part (3), suppose that $f \colon \cpinfty \to X$ is a map.  Write $f_p \colon \cpinfty \to \hpinfty_{(p)}$ for the component map of $f$ corresponding to the prime $p$.  That is, $f_p$ is the composition 
\begin{equation}
\cpinfty \buildrel f \over \to X \to \hpinfty_{(p)}
\end{equation}
where the second map is the natural map arising from the construction of $X$.  Then for any prime $p$ we have the equality
\begin{equation}
\label{eq1:maps}
\dg(f) ~= ~ n_p\, \dg(f_p).
\end{equation}
Since each $n_p$ divides $\dg(f)$, so does their least common multiple $T_X$.  Moreover, by writing $(i_X)_p$ for the component map of $i_X$ corresponding to the prime $p$, \eqref{eq1:maps} implies that for any prime $p$ we have the equalities
\begin{equation}
\label{eq2:maps}
\frac{\dg(f_p)}{\dg(i_X)_p} 
~=~ \frac{\dg(f)/n_p}{T_X/n_p} 
~=~ \frac{\dg(f)}{T_X}.
\end{equation}  
Since there are self-maps of $\hpinfty_{(q)}$ ($q$ any prime) and $\hpinfty_{(0)}$ of degree $\dg(f)/T_X$, one can construct a self-map $g$ of $X$ such that $\dg(g)$ is equal to $\dg(f)/T_X$ and that $f$ is homotopic to $g \circ i_X$.  This proves part (3).

The proof of Theorem \ref{thm:maps} is complete.

\pagebreak

\vspace{1cm}
\noindent
Department of Mathematics \\
University of Illinois at Urbana-Champaign \\
1409 W.\ Green Street \\
Urbana, IL 61801 \\
USA

\vspace{.1cm}
\noindent
\texttt{dyau@math.uiuc.edu}


\begin{thebibliography}{99}

\bibitem{adams}J. F. Adams, \emph{Infinite loop spaces}, Annals of Math. Studies {\bf 90}, Princeton Univ. Press, 1978.

\bibitem{atiyah}M. F. Atiyah, \emph{Power operations in K-theory}, Quart. J. Math. Oxford. {\bf 17} (1966), 165-193.

\bibitem{bk}A. K. Bousfield and D. M. Kan, \emph{Homotopy limits, completions and localizations}, Springer Lecture Notes in Math. {\bf 304}, 1972.

\bibitem{dl}F.-X. Dehon and J. Lannes, 
\emph{Sur les espaces fonctionnels dont la source est le classifiant d'un groupe de Lie compact commutatif}, 
Inst. Hautes Etudes Sci. Publ. Math. {\bf 89} (1999), 127-177. 

\bibitem{imn}K. Ishiguro, J. M$\o$ller, and D. Notbohm, \emph{Rational self-equivalences of spaces in the genus of a product of quaternionic projective spaces}, J. Math. Soc. Japan {\bf 51} (1999), 45-61.

\bibitem{mcgibbon}C. A. McGibbon, \emph{Which group structures on $S^3$ have a maximal torus?}, Springer Lecture Notes in Math. {\bf 657} (1978), 353-360.

\bibitem{nathanson}M. B. Nathanson, \emph{Elementary Methods in Number Theory}, Graduate Texts in Math., vol. 195, Springer-Verlag, New York, 2000.

\bibitem{ns}D. Notbohm and L. Smith, \emph{Fake Lie groups and maximal tori I}, Math. Ann. {\bf 288} (1990), 637-661.

\bibitem{rector}D. L. Rector, \emph{Loop structures on the homotopy type of $S^3$}, Springer Lecture Notes in Math. {\bf 249} (1971), 99-105.

\bibitem{sullivan}D. Sullivan, \emph{Genetics of homotopy theory and the Adams conjecture}, Ann. Math. {\bf 100} (1974), 1-79.

\bibitem{wilkerson}C. Wilkerson, \emph{Classification of spaces of the same n-type for all n}, Proc. Amer. Math. Soc. {\bf 60} (1976), 279-285.

\bibitem{adic}D. Yau, \emph{On adic genus, Postnikov conjugates, and lambda-rings}, \texttt{arXiv:math.AT/0105194}.

\end{thebibliography}
\end{document}